%%%%%%%%%%%% MACROS FOR PAPERS AND REPORTS %%%%%%%%%%%%
%My little macros

\def\ignore#1{}
 
%%%%%%%%%%%%%%%%%%%%%

\newcount\sectnum
\newcount\subsectnum
\newcount\eqnumber

\global\eqnumber=1\sectnum=0

% Equation labels

\def\lab{(\the\sectnum.\the\eqnumber)}

%Example of use: suppose we want to give a label \lgh to an equation
% $$ ......  \xdef\lgh{\lab} \eqnum \show{lgh}$$
% Later refer to Eq. \lgh\ ...
% Note the \ after \lgh; it seems to be needed if we want the equation number
% to be followed by a space; not needed if followed by . or ,

%The next macro is used to display labels in drafts, so that you do
%not have to remember them

\def\show#1{#1}

%The next macro is to be used for final drafts that do not display labels
%\def\show#1{}

%%%%%%%%%%%%%%%%%%%%%

\def\smskip{\vskip 5 pt}
\def\medskip{\vskip 10 pt}
\def\bigskip{\vskip 15 pt}
\def\pn{\par\noindent}
\def\br{\break}

\def\bl{\bigl} 
\def\br{\bigr} 
\def\lf{\left}
\def\ri{\right}

\def\tends{\rightarrow}

\def\frac#1#2{{#1\over #2}}

\def\ol#1{\overline{#1}}

\def\a{\alpha}

\def\m{\mu}
\def\p{\pi}

\def\re{\Re}
\def\rn{\Re^n}

 %break line; horizontal space
\def\tl{\tilde}

\def\old#1{}% invalidates text in braces 
\def\leaderfill{\leaders\hbox to 1em{\hss.\hss}\hfill}
%Example of use: \line{1. Optimality Conditions\leaderfill p.\ 2}

% John's macros

\parindent=2pc
\baselineskip=15pt
\vsize=8.7 true in
\voffset=0.125 true in
\parskip=3pt

% vector/matrix macros

%eqalign macros
\def\minprob#1#2#3{$$\eqalign{&\hbox{minimize\ \ }#1\cr &\hbox{subject to\ \
}#2\cr}\ifnum 0=#3{}\else\eqno(#3)\fi$$}        
     
\def\maxprob#1#2#3{$$\eqalign{&\hbox{maximize\ \ }#1\cr &\hbox{subject to\ \
}#2\cr}\ifnum 0=#3{}\else\eqno(#3)\fi$$}        
     
\def\aligntwo#1#2#3#4#5{$$\eqalign{#1&#2\cr #3&#4\cr}
\ifnum 0=#5{}\else\eqno(#5)\fi$$}
\def\alignthree#1#2#3#4#5#6#7{$$\eqalign{#1&#2\cr #3&#4\cr #5&#6\cr}
\ifnum 0=#7{}\else\eqno(#7)\fi$$}

% Macros to automatically advance equation and other numbers

\def\eqnum{\eqno{\hbox{(\the\sectnum.\the\eqnumber)}\global\advance\eqnumber
by1}}

\def\eqnu{\eqno{\hbox{(\the\sectnum.\the\eqnumber)}\global\advance\eqnumber
by1}}

\newcount\examplnumber
\def\examplnum{\global\advance\examplnumber by1}

\newcount\figrnumber
\def\figrnum{\global\advance\figrnumber by1}

\newcount\propnumber
\def\propnum{\global\advance\propnumber by1}

\newcount\defnumber
\def\defnum{\global\advance\defnumber by1}

\newcount\lemmanumber
\def\lemmanum{\global\advance\lemmanumber by1}

\newcount\assumptionnumber
\def\assumptionnum{\global\advance\assumptionnumber by1}

\newcount\conditionnumber
\def\conditionnum{\global\advance\conditionnumber by1}

\def\exampl{\the\sectnum.\the\examplnumber}
\def\figr{\the\sectnum.\the\figrnumber}
\def\propn{\the\sectnum.\the\propnumber}
\def\defn{\the\sectnum.\the\defnumber}
\def\lemman{\the\sectnum.\the\lemmanumber}
\def\assumptionn{\the\sectnum.\the\assumptionnumber}
\def\condn{\the\sectnum.\the\conditionnumber}

\def\section#1{\goodbreak\vskip 3pc plus 6pt minus 3pt\leftskip=-2pc
   \global\advance\sectnum by 1\eqnumber=1
\global\examplnumber=1\figrnumber=1\propnumber=1\defnumber=1\lemmanumber=1\assumptionnumber=1 \conditionnumber =1%
   \line{\hfuzz=1pc{\hbox to 3pc{\bf %\the\sectnum.\quad
   \vtop{\hfuzz=1pc\hsize=38pc\hyphenpenalty=10000\noindent\uppercase{\the\sectnum.\quad #1}}\hss}}
			\hfill}
			\leftskip=0pc\nobreak\tenf
			\vskip 1pc plus 4pt minus 2pt\noindent\ignorespaces}

% ETP Macros

%\def\section#1{\goodbreak\vskip 3pc plus 6pt minus 3pt\leftskip=-2pc
%   \global\advance\sectnum by 1\eqnumber=1
%   \line{\hfuzz=1pc{\hbox to 3pc{\bf %\the\sectnum.\quad
%   \vtop{\hfuzz=1pc\hsize=38pc\hyphenpenalty=10000\noindent\uppercase{#1}}\hss}}
%                        \hfill}
%                        \leftskip=0pc\nobreak\tenf
%                        \vskip 1pc plus 4pt minus 2pt\noindent\ignorespaces}

\def\sect#1{\noindent\leftskip=-2pc\tenf
   \goodbreak\vskip 1pc plus 4pt minus 2pt
                \global\advance\subsectnum by 1\eqnumber=1
   \line{\hfuzz=1pc{\hbox to 3pc{\bf %\the\sectnum.\quad
   \vtop{\hfuzz=1pc\hsize=38pc\hyphenpenalty=10000\noindent\uppercase{{\bf #1}}}\hss}}
                        \hfill}
   \leftskip=0pc\nobreak\tenf
                        \vskip 1pc plus 4pt minus 2pt\nobreak\noindent\ignorespaces}

\def\subsection#1{\noindent\leftskip=0pc\tenf
   \goodbreak\vskip 1pc plus 4pt minus 2pt
%               \global\advance\subsectnum by 1
   \line{\hfuzz=1pc{\hbox to 3pc{\bf %\the\sectnum.\quad
   \vtop{\hfuzz=1pc\hsize=38pc\hyphenpenalty=10000\noindent{\bf #1}}\hss}}
                        \hfill}
   \leftskip=0pc\nobreak\tenf
                        \vskip 1pc plus 4pt minus 2pt\nobreak\noindent\ignorespaces}
\def\subsubsection#1{\goodbreak\vskip 1pc plus 4pt minus 2pt
   \hfuzz=3pc\leftskip=0pc\noindent\tenit #1 \nobreak\tenf\vskip 6pt plus 1pt
                                minus 1pt\nobreak\ignorespaces\leftskip=0pc}
%
%\def\rthl{Sec. \the\chapnum.\the\sectnum}                      
%\def\rthc{#1}\nobreak\noindent\ignorespaces
%\newcount\sectnum \sectnum=0
%\newcount\subsectnum \subsectnum=0

\def\beginexample#1{\noindent\goodbreak\vskip 6pt plus 1pt minus 1pt
\noindent
  \hbox {\bf Example #1\hss}%\break%\noindent
  \nobreak\vskip 4pt plus 1pt minus 1pt \nobreak\noindent\ninef
  \global\advance
                \leftskip by\parindent\pn}
\def\endexample{\vskip 12pt\tenf\par
  \global\advance\leftskip by -\parindent
  }

\def\beginexercise#1{\noindent\goodbreak\vskip 6pt plus 1pt minus 1pt \noindent\global\normalbaselineskip=12pt
  \hbox {\bf Exercise #1\hss}%\break%\noindent
  \nobreak\vskip 4pt plus 1pt minus 1pt 
  \nobreak\noindent\ninef\global\advance\leftskip
                        by\parindent\pn}
\def\endexercise{\vskip 12pt\tenf\par
  \global\advance\leftskip by -\parindent
  }

\def\beginsection#1{\noindent\goodbreak\vskip 6pt plus 1pt minus 1pt \noindent\global\normalbaselineskip=12pt
  \hbox {\it #1\hss}
  \vskip 0.1pt plus 1pt minus 1pt \nobreak\noindent\ninef\global\advance
                \leftskip by\parindent\noindent\pn}
\def\endsection{\vskip 12pt\tenf\par
  \global\advance\leftskip by -\parindent
}

%

% Header/Title macros

\def\proposition#1{\smskip\pn{\bf Proposition #1}\quad}
\def\proof{\smskip\pn{\bf Proof:}\quad} 
\def\definition#1{\smskip\pn{\bf Definition #1}\quad}

 \def\qed{\quad{\bf
Q.E.D.} \par\bigskip}
\def\ref{\smskip\pn}

\def\chapter#1#2{{\bf \centerline{\helbigbig
{#1}}}\bigskip\bigskip{\bf \centerline{\helbigbig
{#2}}}\bigskip\bigskip} % ex. \chapter{Chapter 1}{Title of chapter}

 % ex. \longchapter{Chapter 1}{Title of chapter}{Title of
 %chapter}

 % ex. \papertitle{Title of paper}{Names of Authors}

\def\longpapertitle#1#2#3{{\bf \centerline{\helbigb
{#1}}}\bigskip{\bf \centerline{\helbigb
{#2}}}\bigskip\bigskip{\centerline{
by}}\bigskip{\bf \centerline{
{#3}}}\bigskip\bigskip} 
% ex. \longpapertitle{First part of title of paper}
%{2nd part of title of paper}{Names of Authors}

% List macros

\def\nitem#1{\smskip\item{#1}}

\newcount\alphanum
\newcount\romnum

\def\alphaenumerate{\ifcase\alphanum \or (a)\or (b)\or (c)\or (d)\or (e)\or
(f)\or (g)\or (h)\or (i)\or (j)\or (k)\fi}
\def\romenumerate{\ifcase\romnum \or (i)\or (ii)\or (iii)\or (iv)\or (v)\or
(vi)\or (vii)\or (viii)\or (ix)\or (x)\or (xi)\fi}

\def\alist{\begingroup\vskip10pt\alphanum=1% alphabetical list
\parskip=2pt\parindent=0pt \leftskip=3pc
\everypar{\llap{{\rm\alphaenumerate\hskip1em}}\advance\alphanum by1}}

\def\nolist{\begingroup\vskip10pt\alphanum=0% numerical list
\parskip=2pt\parindent=0pt \leftskip=3pc
\everypar{\llap{\global\advance\alphanum by1(\the\alphanum)\hskip1em}}}

\def\romlist{\begingroup\vskip10pt\romnum=1% roman list
\parskip=2pt\parindent=0pt \leftskip=5pc
\everypar{\llap{{\rm\romenumerate\hskip1em}}\advance\romnum by1}}

% romlist indents more than alist or nolist and can be used inside them

%Figure, table, and box macros

\long\def\fig#1#2#3{\vbox{\vskip1pc\vskip#1
\prevdepth=12pt \baselineskip=12pt
\vskip1pc
\hbox to\hsize{\hfill\vtop{\hsize=25pc\noindent{\eightbf Figure #2\ }
{\eightpoint#3}}\hfill}}}%Figure space definition. Example of use:
%\fig{16pc}{1.1}{A network with one central processor and a separate
%communication link to each device.}

\long\def\widefig#1#2#3{\vbox{\vskip1pc\vskip#1
\prevdepth=12pt \baselineskip=12pt
\vskip1pc
\hbox to\hsize{\hfill\vtop{\hsize=28pc\noindent{\eightbf Figure #2\ }
{\eightpoint#3}}\hfill}}}

\long\def\table#1#2{\vbox{\vskip0.5pc
\prevdepth=12pt \baselineskip=12pt
\hbox to\hsize{\hfill\vtop{\hsize=25pc\noindent{\eightbf Table #1\ }
{\eightpoint#2}}\hfill}}}

%Running Head Macros
 
\def\rightheadline#1{\headline{\tenrm\hfil #1}}

% Concept Macros

\long\def\leftfig#1#2{\vbox{\smskip\hsize=220pt
\vtop{{\noindent {\bf #1}}}
\smskip
\noindent
\vbox{{\noindent #2}}
}}

\long\def\rightfig#1#2#3{\vbox{\smskip\vskip#1
\prevdepth=12pt \baselineskip=12pt
\hsize=210pt
\smskip
\vbox{\noindent{\eightbold #2}
\hskip1em{\eightpoint#3}}
}}

\long\def\concept#1#2#3#4#5{\bigskip\hrule
\vbox{\hbox{\leftfig{#1}{#2} \hskip3em
\rightfig{#3}{#4}{#5}} \smskip}
\hrule\bigskip}

% Example of Use: \concept{Title of Concept}{Text}
% {Figure size}{Figure number?}{Figure caption}

\long\def\bconcept#1#2#3#4#5#6#7{
\vbox{
\hbox to \hsize{\vtop{\par #1}}
\concept{#2}{#3}{#4}{#5}{#6}
\hbox to \hsize{\vtop{\par #7}}
\smskip}
}

% Example of Use: \bconcept{Preceding text}{Title of Concept}{Text}
% {Figure size}{Figure number}{Figure caption}{Following text}

% same as concept but without the \hrule's; ready to be boxed

% Put inside a box

\def\boxit#1{\vbox{\hrule\hbox{\vrule\kern3pt
                                \vbox{\kern3pt#1\kern3pt}\kern3pt\vrule}\hrule}}
% example of use: \setbox0=\vbox{.... }; \boxit{\box0}
\def\centerboxit#1{$$\vbox{\hrule\hbox{\vrule\kern3pt
                                \vbox{\kern3pt#1\kern3pt}\kern3pt\vrule}\hrule}$$}
% example of use: \setbox0=\vbox{.... }; \centerboxit{\box0}

\long\def\boxtext#1#2{$$\boxit{\vbox{\hsize #1\noindent\strut #2\strut}}$$}
% example of use: \boxtext{462pt}{This is the boxed text.}; 462pt is max length

% Picture macros and examples
%
% figures must be pasted from mcdraw
%
% Look in the 'Windows' menu for the pictures window
% It's like the Scrapbook -- cut and paste pictures
%

\def\picture #1 by #2 (#3){
  \vbox to #2{
    \hrule width #1 height 0pt depth 0pt
    \vfill
    \special{picture #3} % this is the low-level interface
    }
  }
% The first dimension of the picture macro is the width the second is depth

\def\scaledpicture #1 by #2 (#3 scaled #4){{
  \dimen0=#1 \dimen1=#2
  \divide\dimen0 by 1000 \multiply\dimen0 by #4
  \divide\dimen1 by 1000 \multiply\dimen1 by #4
  \picture \dimen0 by \dimen1 (#3 scaled #4)}
  }

%
% Note that you can also say, e.g.,
%  \special{postscript xxx yyy zzz}
% to include PostScript graphics in your documents
%
%Examples of use
%\def\stripes{\picture 2.29in by 1.75in (AWstripes)}
% By executing \stripes 
%\def\annie{\scaledpicture 102pt by 239pt (annie scaled 2000)}
%\def\finder{\picture 260pt by 165pt (screen0 scaled 500)}
%\def\icon{\picture 7in by 7in (icon)}
%Example of use
%\annie
%Example of centered picture \line{\hfil\annie\hfil}

%Figure w/  caption macro
\long\def\captfig#1#2#3#4#5{\vbox{\vskip1pc
\hbox to\hsize{\hfill{\picture #1 by #2 (#3)}\hfill}
\prevdepth=9pt \baselineskip=9pt
\vskip1pc
\hbox to\hsize{\hfill\vtop{\hsize=24pc\noindent{\eightbold Figure #4}
\hskip1em{\eightpoint#5}}\hfill}}}

%Examples of use of Figure macros
%\captfig{8.53pc}{19.9pc}{picturename}{5}{caption.}
%\captfig{2.23in}{2in}{picturename scaled 500}{16}{Caption.}
%The macro centers the picture.
%The first two numbers should be the true width
% and height after the picture has been scaled.
% So if the picture is scaled by 50% (500), the width and height in
% the macro should onw half of what they would be if the picture
% is not scaled (1000).
%
%
%
% Postcript macros

\def\illustration #1 by #2 (#3){
  \vskip#2\hskip#1\special{illustration #3} % this is the low-level interface
    }

\def\scaledillustration #1 by #2 (#3 scaled #4){{
  \dimen0=#1 \dimen1=#2
  \divide\dimen0 by 1000 \multiply\dimen0 by #4
  \divide\dimen1 by 1000 \multiply\dimen1 by #4
  \illustration \dimen0 by \dimen1 (#3 scaled #4)}
  }

% SHADEBOX.BSR MACROS
% Author: Leo@vaxc.cc.monash.edu.au
% Original Source:  Posted by Jimm Herreron <HERRERON@SMCVAX.BITNET>
% Modified from the file SHADEBOX.TEX on 9/30/93 by Becky Kaluza of Blue Sky
% Research to work with Textures 1.5 or later.

\newbox\graybox
\newdimen\xgrayspace
\newdimen\ygrayspace
%
% This macro can be used to typeset some text in a framed box with a
% shaded background. A set of examples can be found at the end of this
% file.
%
% This is a plain \TeX\ file modified for use on the Macintosh with Textures
% 1.5 or later.
%
% The characteristics of the shaded boxes are controlled by the following
% parameters
%
%   \xgrayspace = the space added before and after the text
%   \ygrayspace = the space above and below the text
%   \grayshade  = the gray colour 0 = black 1 = white
%   \linewidth  = the thickness of the border in points
%   \theradius  = the radius of the rounded corners in points
%   \thevskip   = extra \vskip added above and below the shaded box
%                 (applies only to \parashade)
%
%----------------------------------------------------------------------------
%
% The following \TeX code was adapted from previous work by
%
%            Je'ro^me Maillot, maillot@bora.inria.fr
%----------------------------------------------------------------------------
%
% Use the following for one or more words within a line.
%

\def\Textshade#1#2#3#4#5#6{%
    \xgrayspace=#4pt%
    \ygrayspace=#4pt%
    \def\grayshade{#3}%
    \def\linewidth{#5}%
    \def\theradius{#6}%
    \setbox\graybox=\hbox{\surroundboxa{#2}}%
    \hbox{%
    \hbox to 0pt{%
%!    \special{"gsave newpath 0 0 moveto                                %
    \PScommands
}%
    \box\graybox}}%
%
% Use the following for paragraphs.
%
\long%

\long%
\def\Parashade#1#2#3#4#5#6#7{%
    \xgrayspace=#4pt%
    \ygrayspace=#4pt%
    \def\grayshade{#3}%
    \def\linewidth{#5}%
    \def\theradius{#6}%
    \def\thevskip{#7pt}%
    \setbox\graybox=\hbox{\surroundboxb{#2}}%
    \vskip\thevskip%
    \hbox{%
    \hbox to 0pt{%
%!    \special{"gsave newpath 0 0 moveto                                %
    \PScommands
}%
     \box\graybox}%
     \vskip\thevskip%
}%
%----------------------------------------------------------------------------
%
% A pair of box macros. Each builds a slightly oversized box
% containing the text. The text is centred both in the vertical
% horizontal directions.
%
% Use the following for one or more words within a line.
%
\long\def\surroundboxa#1{\leavevmode\hbox{\vtop{%
\vbox{\kern\ygrayspace%
\hbox{\kern\xgrayspace#1%
      \kern\xgrayspace}}\kern\ygrayspace}}}
%
% Use the following for a paragraphs.
%
\long\def\surroundboxb#1{\leavevmode\hbox{\vtop{%
\vbox{\kern\ygrayspace%
\hbox{\kern\xgrayspace\vbox{\advance\hsize-2\xgrayspace#1}%
      \kern\xgrayspace}}\kern\ygrayspace}}}
%----------------------------------------------------------------------------
%
% Here are some simple PostScript routines.
%
% The TeX command \PScommands must be called before any of the
% shading routines can be used.
%
%!\long\def\PScommands{\special{! TeXDict begin
\long\def\PScommands{%
\special{rawpostscript
/sharpbox{%
           newpath
           xmin ymin moveto
           xmin ymax lineto
           xmax ymax lineto
           xmax ymin lineto
           xmin ymin lineto
           closepath 
          }bind def
}%
\special{rawpostscript
/sharpboxnb{%
           newpath
           xmin ymin moveto
           xmin ymax lineto
           xmax ymax lineto
           xmax ymin lineto
%           xmin ymin lineto
%           closepath 
          }bind def
}%
\special{rawpostscript
/sharpboxnt{%
           newpath
           xmin ymax moveto
           xmin ymin lineto
           xmax ymin lineto
           xmax ymax lineto
%           xmin ymin lineto
%           closepath 
          }bind def
}%
\special{rawpostscript
/roundbox{%
           newpath
           xmin radius add ymin moveto
           xmax ymin xmax ymax radius arcto
           xmax ymax xmin ymax radius arcto
           xmin ymax xmin ymin radius arcto
           xmin ymin xmax ymin radius arcto 16 {pop} repeat
           closepath
          }bind def
}%
\special{rawpostscript
/sharpcorners{%
               sharpbox gsave grayshade setgray fill grestore 
               linewidth setlinewidth stroke
              }bind def
}%
\special{rawpostscript
/sharpcornersnt{%
               sharpboxnt gsave grayshade setgray fill grestore 
               linewidth setlinewidth stroke
              }bind def
}%
\special{rawpostscript
/sharpcornersnb{%
               sharpboxnb gsave grayshade setgray fill grestore 
               linewidth setlinewidth stroke
              }bind def
}%
\special{rawpostscript
/roundcorners{%
               roundbox gsave grayshade setgray fill grestore 
               linewidth setlinewidth stroke
              }bind def
}%
\special{rawpostscript
/plainbox{%
           sharpbox grayshade setgray fill 
          }bind def
}%
% Here are the two new options
%
\special{rawpostscript
/roundnoframe{%
               roundbox grayshade setgray fill 
              }bind def
}%
\special{rawpostscript
/sharpnoframe{%
               sharpbox grayshade setgray fill 
              }bind def
}%
%!end}%
}%
%
% The \PScommands macro must be invoked before the shaded box macros.
%
%!\PScommands
% To use this, type \textshade{plainbox} or \textshade{roundbox} or
% \textshade{sharpbox}

\def\pshade#1{\Parashade{sharpcorners}{#1}{0.95}{10}{0.5}{10}{10}}

%%%%% BOXES FOR TEXSHOP %%%%%

\def\boxit#1{\vbox{\hrule\hbox{\vrule\kern3pt
                                \vbox{\kern3pt#1\kern3pt}\kern3pt\vrule}\hrule}}
% example of use: \setbox0=\vbox{.... }; \boxit{\box0}

\def\boxitnb#1{\vbox{\hrule\hbox{\vrule\kern3pt
                                \vbox{\kern3pt#1\kern3pt}\kern3pt\vrule}}}

\def\boxitnt#1{\vbox{\hbox{\vrule\kern3pt
                                \vbox{\kern3pt#1\kern3pt}\kern3pt\vrule}\hrule}}

\long\def\boxtext#1#2{$$\boxit{\vbox{\hsize #1\noindent\strut #2\strut}}$$}
% example of use: \boxtext{462pt}{This is the boxed text.}; 462pt is max length

% example of use: \boxtext{462pt}{This is the boxed text.}; 462pt is max length

% example of use: \boxtext{462pt}{This is the boxed text.}; 462pt is max length

\def\texshopbox#1{\boxtext{462pt}{\vskip-1.5pc\pshade{\vskip-1.0pc#1\vskip-2.0pc}}}

%***************************************************
%         FONTS
%***************************************************

% ROMAN
%
%
%
%
%
%
%
%
\font\helbigbig=cmr10 scaled 2500%
\font\helbigb=cmbx10 scaled 1500%
\font\eightbold=cmbx8%

\def\tenf{\hel}%
\def\tenit{\heli}%
\def\ninef{\ninehel}%
\def\nineit{\nineheli}%
%
%

%  FONT FAMILIES

\font\tenrm=cmr10%
\font\teni=cmmi10%
\font\tensy=cmsy10%
\font\tenbf=cmbx10%
\font\tentt=cmtt10%
\font\tenit=cmti10%
\font\tensl=cmsl10%

\def\tenpoint{\def\rm{\fam0\tenrm}%
\textfont0=\tenrm%
\textfont1=\teni%
\textfont2=\tensy%
\textfont\itfam=\tenit%
\textfont\slfam=\tensl%
\textfont\ttfam=\tentt%
\textfont\bffam=\tenbf%
\scriptfont0=\sevenrm%
\scriptfont1=\seveni%
\scriptfont2=\sevensy%
%\scriptfont3=\tenex%
\scriptscriptfont0=\sixrm%
\scriptscriptfont1=\sixi%
\scriptscriptfont2=\sixsy%
%\scriptscriptfont3=\tenex%
\def\it{\fam\itfam\tenit}%
\def\tt{\fam\ttfam\tentt}%
\def\sl{\fam\slfam\tensl}%
\scriptfont\bffam=\sevenbf%
\scriptscriptfont\bffam=\sixbf%
\def\bf{\fam\bffam\tenbf}%
\normalbaselineskip=18pt%
\normalbaselines\rm}%

\font\ninerm=cmr9%
\font\ninebf=cmbx9%
\font\nineit=cmti9%
\font\ninesy=cmsy9%
\font\ninei=cmmi9%
\font\ninett=cmtt9%
\font\ninesl=cmsl9%

\def\ninepoint{\def\rm{\fam0\ninerm}%
\textfont0=\ninerm%
\textfont1=\ninei%
\textfont2=\ninesy%
\textfont\itfam=\nineit%
\textfont\slfam=\ninesl%
\textfont\ttfam=\ninett%
\textfont\bffam=\ninebf%
\scriptfont0=\sixrm%
\scriptfont1=\sixi%
\scriptfont2=\sixsy%
%\scriptfont3=\tenex%
\def\it{\fam\itfam\nineit}%
\def\tt{\fam\ttfam\ninett}%
\def\sl{\fam\slfam\ninesl}%
\scriptfont\bffam=\sixbf%
\scriptscriptfont\bffam=\fivebf%
\def\bf{\fam\bffam\ninebf}%
\normalbaselineskip=16pt%
\normalbaselines\rm}%

\font\eightrm=cmr8%
\font\eighti=cmmi8%
\font\eightsy=cmsy8%
\font\eightbf=cmbx8%
\font\eighttt=cmtt8%
\font\eightit=cmti8%
\font\eightsl=cmsl8%

\def\eightpoint{\def\rm{\fam0\eightrm}%
\textfont0=\eightrm%
\textfont1=\eighti%
\textfont2=\eightsy%
\textfont\itfam=\eightit%
\textfont\slfam=\eightsl%
\textfont\ttfam=\eighttt%
\textfont\bffam=\eightbf%
\scriptfont0=\sixrm%
\scriptfont1=\sixi%
\scriptfont2=\sixsy%
%\scriptfont3=\tenex%
\scriptscriptfont0=\fiverm%
\scriptscriptfont1=\fivei%
\scriptscriptfont2=\fivesy%
%\scriptscriptfont3=\tenex%
\def\it{\fam\itfam\eightit}%
\def\tt{\fam\ttfam\eighttt}%
\def\sl{\fam\slfam\eightsl}%
%\scriptfont\bffam=\sixbf%
\scriptscriptfont\bffam=\fivebf%
\def\bf{\fam\bffam\eightbf}%
\normalbaselineskip=14pt%
\normalbaselines\rm}%

\font\sevenrm=cmr7%
\font\seveni=cmmi7%
\font\sevensy=cmsy7%
\font\sevenbf=cmbx7%

\def\sevenpoint{%
   \def\rm{\sevenrm}\def\bf{\sevenbf}%
   \def\smc{\sevensmc}\baselineskip=12pt\rm}%

\font\sixrm=cmr6%
\font\sixi=cmmi6%
\font\sixsy=cmsy6%
\font\sixbf=cmbx6%

\fontdimen13\tensy=2.6pt%
\fontdimen14\tensy=2.6pt%
\fontdimen15\tensy=2.6pt%
\fontdimen16\tensy=1.2pt%
\fontdimen17\tensy=1.2pt%
\fontdimen18\tensy=1.2pt%       

\def\tenf{\tenpoint}%
\def\ninef{\ninepoint}%
%

%%%%%%%%%%%% END OF MACROS %%%%%%%%%%%%

%\input TEX SHOP_small_baseline.tex

\def\section#1{\goodbreak\vskip 3pc plus 6pt minus 3pt\leftskip=-2pc
   \global\advance\sectnum by 1\eqnumber=1\subsectnum=0%
\global\examplnumber=1\figrnumber=1\propnumber=1\defnumber=1\lemmanumber=1\assumptionnumber=1 \conditionnumber =1%
   \line{\hfuzz=1pc{\hbox to 3pc{\bf %\the\sectnum.\quad
   \vtop{\hfuzz=1pc\hsize=38pc\hyphenpenalty=10000\noindent\uppercase{\the\sectnum.\quad #1}}\hss}}
			\hfill}
			\leftskip=0pc\nobreak\tenf
			\vskip 1pc plus 4pt minus 2pt\noindent\ignorespaces}
\def\subsection#1{\noindent\leftskip=0pc\tenf
   \goodbreak\vskip 1pc plus 4pt minus 2pt
               \global\advance\subsectnum by 1
   \line{\hfuzz=1pc{\hbox to 3pc{\bf \the\sectnum.\the\subsectnum.
   \vtop{\hfuzz=1pc\hsize=38pc\hyphenpenalty=10000\noindent{\bf #1}}\hss}}
                        \hfill}
   \leftskip=0pc\nobreak\tenf
                        \vskip 1pc plus 4pt minus 2pt\nobreak\noindent\ignorespaces}

%\newcount\conditionnumber
%\def\conditionnum{\global\advance\conditionnumber by1}
%\def\conditionn{\the\sectnum.\the \conditionnumber} 

%%%%%%%%%% REDEFINITION OF BOX SPACING %%%%%%%%%%%%%%%%

\def\texshopbox#1{\boxtext{462pt}{\vskip-1.5pc\pshade{\vskip-1.0pc#1\vskip-2.0pc}}}

%%%%%%%%%%%%%%%%%%%%%%%%%%%%%%%%%%%%%%%%%%%%

\long\def\fig#1#2#3{\vbox{\vskip1pc\vskip#1
\prevdepth=12pt \baselineskip=12pt
\vskip1pc
\hbox to\hsize{\hfill\vtop{\hsize=30pc\noindent{\eightbf Figure #2\ }
{\eightpoint#3}}\hfill}}}

\def\show#1{}

\rightheadline{\botmark}

\pageno=1

\def\longpapertitle#1#2#3{{\bf \centerline{\helbigb
{#1}}}\medskip{\bf \centeline{\helbigb
{#2}}}\bigskip{\bf \centerline{
{#3}}}\bigskip}

\vskip-3pc

\def\xstar{X^{\raise0.04pt\hbox{\sevenpoint *}} }

\def\jstar{J^{\raise0.04pt\hbox{\sevenpoint *}} }
\def\qstar{Q^{\raise0.04pt\hbox{\sevenpoint *}} }

\rightheadline{\botmark}

\pageno=1

\rightheadline{\botmark}

\pageno=1

\rightheadline{\botmark}

\pn {\bf May 2021}%\hfill{\bf  DRAFT NOTE}%
\bigskip \bigskip\bigskip

\bigskip

\def\longpapertitle#1#2#3{{\bf \centerline{\helbigb
{#1}}}\medskip{\bf \centerline{\helbigb
{#2}}}\bigskip{\bf \centerline{
{#3}}}\bigskip}

\vskip-3pc

\longpapertitle{On-Line Policy Iteration for Infinite }{Horizon Dynamic Programming}{{Dimitri Bertsekas\footnote{\dag}{\ninepoint Fulton Professor of Computational Decision Making, ASU, Tempe, AZ, and McAfee Professor of Engineering, MIT, Cambridge, MA.}}}

%\vskip-1.5pc

\centerline{\bf Abstract}

\smskip
\pn In this paper we  propose an on-line policy iteration (PI) algorithm for finite-state infinite horizon discounted dynamic programming, whereby the policy improvement operation is done on-line, only for the states that are encountered during operation of the system. This allows the continuous updating/improvement of the current policy, thus resulting in a form of on-line PI that incorporates the improved controls into the current policy as new states and controls are generated. The algorithm converges in a finite number of stages to a type of locally optimal policy, and suggests the possibility of variants of PI and multiagent PI where the policy improvement is simplified. Moreover, the algorithm can be used with on-line replanning, and is also well-suited for on-line PI algorithms with value and policy approximations.

\vskip-1pc

\section{Simplified On-Line Policy Iteration for Discounted Problems}
\vskip-0.5pc
\pn We introduce new on-line variants of the classical policy iteration (PI) algorithm for finite-state discounted infinite horizon dynamic programming (DP) problems. The common characteristic of these variants is that, in addition to being suitable for on-line implementation, they are simplified in two ways: 
\nitem{(a)} They perform policy improvement operations only for the states that are encountered during the on-line operation of the system.
\nitem{(b)} The policy improvement operation is simplified in that it uses approximate minimization over the Q-factors of the current policy at the current state.
\smskip
\pn Despite these simplifications, we show that our algorithms generate a sequence of improved policies, which converge to a policy with a local optimality property. Moreover, with an enhancement of the policy improvement operation, which involves a form of exploration, they converge to a globally optimal policy. 

The motivation for our work comes from the rollout algorithm; see the author's reinforcement learning books [Ber19] and [Ber20], which provide many additional references. This algorithm starts from some available ``base policy" and implements on-line an improved ``rollout policy," which is the one that would be obtained from the first iteration of the PI algorithm starting from the base policy. In the algorithm of the present paper, the data accumulated from the rollout implementation is used to improve on-line the base policy, and to asymptotically obtain a policy that is either locally or globally optimal.

We assume a discrete-time
dynamic system with states $1,\ldots,n$, and we use a transition probability notation. We denote states by the symbol $x$ and successor states by the symbol $y$. The control/action is denoted by $u$, and is constrained to take values in a given finite constraint set $U(x)$, which may depend on the current state $x$. The use of a control $u$ at state
$x$ specifies the transition probability
$p_{xy}(u)$ to the next state  $y$, at a cost $g(x,u,y)$.

A policy $\pi=\{\m_0,\m_1,\ldots\}$ is a sequence of functions from state to control that satisfies the control constraint, i.e.,  $\m_k(x)\in U(x)$ for all $x$ and $k$. Given a policy $\p$ and an initial state $x_0$, the  system becomes a Markov chain whose generated trajectory under $\p$, denoted $\{x_0,x_1,\ldots\}$, has a well-defined probability distribution. The corresponding total expected cost is
$$J_\p(x_0) = \lim_{N\tends\infty} E\lf\{\sum_{k=0}^{N-1}
\a^k g\bl(x_k,\m_k(x_k),x_{k+1}\br)\ \Big|\ x_0,\,\p\ri\},\qquad x_0=1,\ldots,n,$$
where $\a<1$ is the discount factor. The above expected value is taken with respect to the joint distribution of the states $x_1,x_2,\ldots$, conditioned on the initial state being $x_0$ and the use of $\p$. The optimal cost starting from a state
$x$, i.e., the minimum of
$J_\p(x)$ over all policies $\p$, is denoted by $\jstar(x)$. We will view $\jstar$ as the vector of the $n$-dimensional space $\rn$ that has components $\jstar(1),\ldots,\jstar(n)$.\footnote{\dag}{\ninepoint  In our notation, $\re^n$ is the $n$-dimensional Euclidean space and all vectors in $\re^n$ are viewed as column vectors. Moreover, all vector inequalities $J\le J'$ are meant be componentwise, i.e., $J(x)\le J'(x)$ for all $x=1,\ldots,n$.}

A stationary policy is a policy of the form $\pi =
\{\mu,\mu,\ldots\}$, and for
brevity, it is referred to as the ``policy $\m$." The cost of $\m$ starting from state $x$ is denoted by $J_\m(x)$, and is given by
$$J_\m(x_0) = \lim_{N\tends\infty} E\lf\{\sum_{k=0}^{N-1}
\a^k g\bl(x_k,\m(x_k),x_{k+1}\br)\ \Big|\ x_0,\,\m\ri\},\qquad x_0=1,\ldots,n.$$
We can view $J_\m$ as the vector in $\rn$ that has components $J_\m(1),\ldots,J_\m(n)$.  We say that
$\m$ is optimal if 
$$J_\m(x)=\jstar(x)=\min_\p J_\p(x),\qquad x=1,\ldots,n.$$
It is well known that there exists an optimal stationary policy; see e.g., the books [Ber12], [Put94], which provide an extensive analysis of discounted finite-state infinite horizon DP problems.

The theory and algorithms for our problem are conveniently stated with the use of abstract notation, as in the book [Ber18]. In particular, for each policy $\m$, we introduce the operator $T_\m:\rn\mapsto\rn$, which maps a vector $J\in\rn$ to the vector $T_\m J\in\rn$ that has components 
$$(T_\m J)(x)=\sum_{y=1}^np_{xy}\big(\m(x)\big)\big(g(x,\m(x),y)+\a J(y)\big),\qquad x=1,\ldots,n.\eqnum\show{oneo}$$
We also introduce the operator $T:\rn\mapsto\rn$ defined by
$$(T J)(x)=\min_{u\in U(x)}\sum_{y=1}^np_{xy}(u)\big(g(x,u,y)+\a J(y)\big),\qquad x=1,\ldots,n.\xdef\discountstar{\lab}\eqnum\show{oneo}$$

An important property is that $T_\m$ and $T$ are monotone, i.e., that for all $J,J'\in\rn$,
$$T_\m J\le T_\m J',\qquad TJ\le TJ',\qquad \hbox{if }J\le J'.$$ 
Another important property is that $T_\m$ and $T$ are sup-norm contractions, so that the costs $J_\m(x)$, $x=1,\ldots,n,$ are the unique solution of Bellman's equation
$$J_\m(x)=\sum_{y=1}^np_{xy}\big(\m(x)\big)\big(g(x,\m(x),y)+\a J_\m(y)\big),\qquad x=1,\ldots,n,\eqnum\show{oneo}$$ 
which can also be written as the fixed point equation $J_\m=T_\m J_\m$. Similarly, the optimal costs $\jstar(x)$, $x=1,\ldots,n$, are the unique solution of Bellman's equation
$$J^*(x)=\min_{u\in U(x)}\sum_{y=1}^np_{xy}(u)\big(g(x,u,y)+\a J^*(y)\big),\qquad x=1,\ldots,n,\eqnum\show{oneo}$$
or $\jstar=T\jstar$. A consequence of this is that the following optimality conditions hold
$$T_\m \jstar=T\jstar\qquad \hbox{if and only if}\qquad \hbox{$\m$ is optimal},\eqnum\show{oneo}$$
$$T_\m J_\m=TJ_\m\qquad \hbox{if and only if}\qquad \hbox{$\m$ is optimal}.\xdef\optcondtwo{\lab}\eqnum\show{oneo}$$
The contraction property also implies that the value iteration (VI) algorithms 
$$J^{k+1}=T_\m J^k,\qquad J^{k+1}=TJ^k$$ 
generate sequences $\{J^k\}$ that converge  to $J_\m$ and $J^*$, respectively, from any starting vector $J^0\in\rn$. 

Policy iteration (PI) is a major alternative to VI, and generates a sequence of policies. 
In the classical form of the algorithm, the current policy  $\m$ is improved by finding $\tl \m$ that satisfies $$T_{\tl \m}J_\m=TJ_\m$$
 [i.e., by minimizing for all $x$ in the right-hand side of Eq.\ \discountstar\ with $J_\m$ in place of $J$]. The improved policy $\tl\m$ is evaluated by solving the linear system of equations $J_{\tl\m}=T_{\tl\m}J_{\tl\m}$, and then $(J_{\tl\m},\tl\m)$ becomes the new cost vector-policy pair, which is used to start a new iteration. Thus the PI algorithm starts with a policy $\m^0$ and generates a sequence $\{\m^k\}$ according to
$$J_{\m^k}=T_{\m^k}J_{\m^k},\qquad T_{\m^{k+1}}J_{\m^k}=TJ_{\m^k},\eqnum\show{oneo}$$
where on the left we have the policy evaluation equation, and on the right we have the policy improvement equation.

The preceding PI algorithm has several weaknesses, which make it unsuitable for problems with a large number of states $n$. The principal difficulty is that at each iteration $k$, the current policy $\m^k$ must be evaluated at all states, i.e., $J_{\m^k}(x)$ must be computed for all $x$. A second difficulty is that the policy improvement operation 
$$\m^{k+1}(x)\in\arg\min_{u\in U(x)}\sum_{y=1}^np_{xy}(u)\big(g(x,u,y)+\a J_{\m^k}(y)\big),\qquad x=1,\ldots,n,\eqnum\show{oneo}$$
must be performed at each state $x$,  which makes it unsuitable for problems where the number of state-control pairs $(x,u)$ is  large. The purpose of this paper is to propose variants of the PI algorithm that alleviate both of these difficulties.

\vskip-1pc

\section{On-Line Variants of Policy Iteration}
  
 \pn We introduce a PI algorithm, which starts at time 0 with a state-policy pair $(x_0,\m^0)$ and generates on-line a sequence of state-policy pairs $(x_k,\m^k)$. We view $x_k$ as the current state of a system operating on line under the influence of the policies $\m^1,\m^2,\ldots$. 
In our algorithm, {\it $\m^{k+1}$ may differ from $\m^k$ only at state $x_k$\/}. The control $\m^{k+1}(x_k)$ and the state $x_{k+1}$ are generated as follows:

We consider the Q-factors
$$Q_{\m^k}(x_k,u)=\sum_{y=1}^np_{x_ky}(u)\big(g(x_k,u,y)+\a J_{\m^k}(y)\big),\xdef\qfactork{\lab}\eqnum\show{oneo}$$
and we select the control $\m^{k+1}(x_k)$ from within the constraint set $U(x_k)$ with sufficient accuracy to satisfy the sequential improvement condition
$$Q_{\m^k}\big(x_k,\m^{k+1}(x_k)\big)\le J_{\m^k}(x_k),\xdef\seqcondsimpl{\lab}\eqnum\show{oneo}$$
with strict inequality whenever this is possible.\footnote{\dag}{\ninepoint By this we mean that if $\min_{u\in U(x_k)}Q_{\m^k}(x_k,u)<J_{\m^k}(x_k)$ we select a control $u_k$ that satisfies
$$Q_{\m^k}(x_k,u_k)<J_{\m^k}(x_k),\xdef\strictcondsimpl{\lab}\eqnum\show{oneo}$$
and set $\m^{k+1}(x_k)=u_k$, and otherwise we set $\m^{k+1}(x_k)=\m^k(x_k)$ [so Eq.\ \seqcondsimpl\ is satisfied]. Such a control selection may be obtained by a number of schemes, including brute force calculation and random search based on Bayesian optimization. The needed values of the Q-factor $Q_{\m^k}$ and cost $J_{\m^k}$ may be obtained in several ways, depending on the problem at hand, including by on-line simulation.}
For $x\ne x_k$ the policy controls are not changed:
$$\m^{k+1}(x)=\m^k(x)\ \hbox{ for all $x\ne x_k$.}$$
 The next state $x_{k+1}$ is generated randomly according to the transition probabilities $p_{x_kx_{k+1}}\big(\m^{k+1}(x_k)\big)$.

We first show that the current policy is monotonically improved.

\xdef\propone{\propn}\propnum\show{myproposition}

\texshopbox{\proposition{\propone:} We have 
$$J_{\m^{k+1}}(x)\le J_{\m^k}(x),\qquad\hbox{for all $x$ and $k$},$$
with strict inequality for $x=x_k$ (and possibly other values of $x$) if $\min_{u\in U(x_k)}Q_{\m^k}(x_k,u)<J_{\m^k}(x_k)$.
}
\proof The policy update is done under the condition  \seqcondsimpl. By using the monotonicity of $T_{\m^{k+1}}$, we have for all $\ell\ge1$,
$$T_{\m^{k+1}}^{\ell+1}J_{\m^k}\le T_{\m^{k+1}}^{\ell}J_{\m^k}\le J_{\m^k},\xdef\monotoneineq{\lab}\eqnum\show{oneo}$$ 
so by taking the limit as $\ell\to\infty$ and by using the convergence property of VI ($T_{\m^{k+1}}^{\ell}J\to J_{\m^{k+1}}$ for any $J$), we obtain $J_{\m^{k+1}}\le J_{\m^k}$. Moreover,  the algorithm selects $\m^{k+1}(x_k)$ so that
$$(T_{\m^{k+1}}J_{\m^{k}})(x_k)=Q_{\m^k}(x_k,u_k)< J_{\m^k}(x_k)\qquad \hbox{if }\ \min_{u\in U(x_k)}Q_{\m^k}(x_k,u)<J_{\m^k}(x_k),$$
[cf.\ Eq.\ \strictcondsimpl],  so that by using Eq.\ \monotoneineq, we have $J_{\m^{k+1}}(x_k)< J_{\m^k}(x_k)$. \qed

\subsubsection{Convergence to a Locally Optimal Policy}

\pn We next discuss the convergence and optimality properties of the algorithm. We introduce a definition of local optimality of a policy, whereby the policy selects controls optimally only within a subset of states.

%\xdef\defone{\propn}\propnum\show{myproposition}

\texshopbox{\definition{2.1:}Given a subset of states $X$ and a policy $\m$, we say that  {\it $\m$ is locally optimal over $X$}  if $\m$ is optimal for the problem where the control is restricted to take the value $\m(x)$ at the states $x\notin X$, and is allowed to take any value $u\in U(x)$ at the states $x\in X$.}

Roughly speaking, $\m$ is locally optimal over $X$, if $\m$ is acting optimally within $X$, but under the (incorrect) assumption that once the state of the system gets to a state $x$ outside $X$, there will be no option to select control other than $\m(x)$. Thus if the choices of $\m$ outside $X$ are poor, its choices within $X$ may also be poor.

Mathematically,  $\m$ is locally optimal over $X$ if
$$J_\m(x)=\min_{u\in U(x)}\sum_{y=1}^np_{xy}(u)\big(g(x,u,y)+\a J_\m(y)\big),\qquad \hbox{for all } x\in X,$$
$$J_\m(x)=\sum_{y=1}^np_{xy}\big(\m(x)\big)\Big(g\big(x,\m(x),y\big)+\a J_\m(y)\Big),\qquad \hbox{for all } x\notin X,$$
which can be written compactly as
$$(T_\m J_\m)(x)=(TJ_\m)(x),\qquad \hbox{for all $x\in X$}.\xdef\localoptimality{\lab}\eqnum\show{oneo}$$
Note that this is different than (global) optimality of $\m$, which holds if and only if the above condition holds for all $x=1,\ldots,n$, rather than just for $x\in X$ [cf.\ Eq.\ \optcondtwo]. However, it can be seen from Definition 2.1 that a (globally) optimal policy is also locally optimal within any subset of states. 

Our main convergence result is the following.

\xdef\proptwo{\propn}\propnum\show{myproposition}

\texshopbox{\proposition{\proptwo:} Let $\ol X$ be the subset of states that are repeated infinitely often within the sequence $\{x_k\}$. Then the corresponding sequence $\{\m^k\}$ converges finitely to some policy $\ol \m$ in the sense that $\m^k=\ol\m$ for all $k$ after some index $\ol k$. Moreover $\ol \m$ is locally optimal within $\ol X$, while $\ol X$ is invariant under $\ol \m$, in the sense that
$$p_{xy}\big(\ol \m(x)\big)=0\qquad\hbox{for all $x\in \ol X$ and $y\notin \ol X$}.$$
}

\proof The cost function sequence $\{J_{\m^k}\}$ is monotonically nondecreasing (cf.\ Prop.\ \propone). The number of policies $\m$ is finite in view of the finiteness of the state and control spaces. Therefore, the number of corresponding functions $J_\m$ is also finite,  so $J_{\m^k}$ converges in a finite number of steps to some $\skew5\bar J$, which in view of the algorithm's construction [selecting $u_k=\m^k(x_k)$  if $\min_{u\in U(x_k)}Q_{\m^k}(x_k,u)=J_{\m^k}(x_k)$; cf.\ Eq.\ \strictcondsimpl], implies that $\m^k$ will remain unchanged at some $\ol \m$ with $J_{\ol\m}=\skew5\bar J$ after some sufficiently large $k$. 

We will show that the local optimality condition \localoptimality\ holds for $X=\ol X$ and $\m=\ol\m$. In particular, we have $x_k\in \ol X$ and $\m^k=\ol\m$ for all $k$ greater than some index, while for every $x\in\ol X$, we have $x_k=x$ for infinitely many $k$. It follows that for all $x\in\ol X$, 
$$Q_{\ol\m}\big(x,\ol\m(x)\big)= J_{\ol\m}(x),\xdef\conditionone{\lab}\eqnum\show{oneo}$$
while by the construction of the algorithm, 
$$Q_{\ol\m}\big(x,u\big)\ge  J_{\ol\m}(x),\qquad \hbox{for all $u\in U(x)$,}\xdef\conditiontwo{\lab}\eqnum\show{oneo}$$
since the reverse would imply that $\m^{k+1}(x)\ne \m^k(x)$ for infinitely many $k$ [cf.\ Eq.\ \strictcondsimpl]. Condition \conditionone\ can be written as $J_{\ol\m}(x)=(T_{\ol\m}J_{\ol\m})(x)$ for all $x\in\ol X$, and combined  with Eq.\ \conditiontwo, implies that $(T_{\ol\m}J_{\ol\m})(x)=(TJ_{\ol\m})(x)$ for all $x\in\ol X$. This is the local optimality condition \localoptimality\ with  $X=\ol X$ and $\m=\ol\m$. 

To show that $\ol X$ is invariant under $\ol\m$, we argue by contradiction: if this were not so, there would exist a state $x\in \ol X$ and a state $y\notin\ol X$ such that
$$p_{xy}\big(\ol\m(x)\big)>0,$$
implying that $y$ would be generated following the occurrence of $x$ infinitely often within the sequence $\{x_k\}$, and hence would have to belong to $\ol X$ (by the definition of $\ol X$). \qed

Note an implication of the invariance property of the set $\ol X$ shown in the preceding proposition. We have that $\ol \m$ is (globally) optimal under the assumption that for every policy there is no strict subset of states that is invariant.

\subsubsection{A Counterexample to Global Optimality}

\pn The following deterministic example\footnote{$\,$\dag}{\ninepoint Thanks are due to Yuchao Li for constructing this example and for commenting on other aspects of the paper.} shows that the policy $\ol \m$ obtained by the algorithm need not be (globally) optimal. 
Here there are three states 1, 2, and 3. From state 1 we can go to state 2 at cost 1, and to state 3 at cost 0, from state 2 we can go to states 1 and 3 at cost 0, and from state 3 we can go to state 2 at cost 0 or stay in 3 at a high cost (say 10). The discount factor is $\a=0.9$. Then it can be verified that the optimal policy is 
$$\m^*(1): \hbox{Go to 3},\qquad \m^*(2): \hbox{Go to 3},\qquad \m^*(3): \hbox{Go to 2},$$
with optimal costs
$$J^*(1)=J^*(2)=J^*(3)=0,$$
 while the policy
$$\ol \m(1): \hbox{Go to 2},\qquad \ol \m(2): \hbox{Go to 1},\qquad \ol \m(3): \hbox{Stay at 3},$$
is strictly suboptimal, but is locally optimal over the set of states $\ol X=\{1,2\}$. Moreover our on-line PI algorithm, starting from state 1 and the policy $\m^0=\ol\m$, oscillates between the states 1 and 2, and leaves the policy $\m^0$ unchanged. Note also that $\ol X$ is invariant under $\ol\m$, consistently with Prop.\ \proptwo.

\subsubsection{On-Line Variants of Policy Iteration with Global Optimality  Properties}

\pn To address the local versus global convergence issue illustrated by the preceding example,  we consider an alternative scheme, whereby in addition to $u_k$, we generate an additional control at a randomly chosen state $\ol x_k\ne x_k$.\footnote{\dag}{\ninepoint It is also possible to choose multiple additional states at time $k$ for a policy improvement operation, and this is well-suited for the use of parallel computation.}
In particular, assume that at each time $k$, in addition to $u_k$ and $x_{k+1}$ that are generated according to Eq.\ \strictcondsimpl, the algorithm generates randomly another state $\ol x_{k}$ (all states are selected with positive probability), performs a policy improvement operation at that state as well, and modifies accordingly $\m^{k+1}(\ol x_{k})$. Thus, in addition to a policy improvement operation at each state within the generated sequence $\{x_k\}$, there is an additional policy improvement operation at each state within the randomly generated sequence $\{\ol x_k\}$. 

Because of the random mechanism of selecting $\ol x_{k}$, it follows that at every state there will be a policy improvement operation infinitely often, which implies that the policy $\ol \m$ ultimately obtained is (globally) optimal. Note also that {\it we may view the random generation of the sequence $\{\ol x_k\}$ as a form of exploration\/}. The probabilistic mechanism for generating the random sequence $\{\ol x_k\}$ may be guided by some heuristic reasoning, which aims to explore states with high cost improvement potential.

\vskip-1pc

\section{Concluding Remarks}
 
\pn We have developed a new on-line PI algorithm, which has a cost improvement property at each iteration, and offers some optimality guarantees. The algorithm is, to our knowledge, the first on-line proposal of the PI algorithm, which in its exact form relies on off-line computation. The algorithm is motivated and inspired by the rollout algorithm (a single iteration of PI), which is well-suited for on-line implementation; see the books  [BeT96], [Ber17], [Ber19], [Ber20], and the references given there. Note that if there were no policy updates at all, our algorithm would become equivalent to the rollout algorithm with base policy $\m^0$. Thus our algorithm may be viewed as an enhanced version of a simplified rollout algorithm, which uses information collected on-line to improve the current base policy (at essentially no additional computational cost). A logical conclusion is that our algorithm should perform no worse than the  rollout algorithm, which is known to perform well in practice, as it can be viewed as a single step of Newton's method, which underlies the PI algorithm [Kle68], [PoA69], [Hew71], [PuB78], [PuB79], [Ber20]. 

The structure of our algorithm suggests interesting possibilities for  exploration, as well as combinations with approximation in policy and value space schemes, whereby policy and cost function approximations are constructed using the data that is collected on-line. Value and policy space approximation could also be applied to problems with infinite state and control spaces. There are many possibilities for the use of parallel computation within these contexts. Moreover, the on-line structure of our algorithm makes it suitable for adaptive control contexts, where some of the system and cost parameters may be changing over time. 

Other possible variations of our algorithm include PI schemes with multistep lookahead, as well as optimistic variants, whereby the Q-factors \qfactork\ are evaluated approximately by using truncated simulation. Extensions to other infinite horizon DP models, such as semi-Markov and average cost problems, are also possible. A particularly interesting class of models is deterministic and stochastic shortest path problems with an additional termination state $t$. For such problems, each generated system trajectory consists of a finite number of states and terminates at $t$, so the variant given at the end of the preceding section, which guarantees global optimality in the limit, does not apply. However, it is possible to construct a globally optimal policy by restarting the system from a randomly chosen initial state each time it reaches $t$.  These possibilities are subjects for further research and computational experimentation. 
 
\vskip-1pc

\section{References}

\ref [BeT96] Bertsekas, D.\ P., and Tsitsiklis, J.\ N., 1996.\ Neuro-Dynamic Programming,
Athena Scientific, Belmont, MA.

\ref[Ber12] Bertsekas, D.\ P., 2012.\ Dynamic Programming and Optimal Control, Vol.\ II, 4th edition, Athena Scientific, Belmont, MA.

\ref[Ber17] Bertsekas, D.\ P., 2017.\ Dynamic Programming and Optimal Control, Vol.\ I, 4th Edition, Athena Scientific, Belmont, MA.

\ref[Ber18] Bertsekas, D.\ P., 2018.\ Abstract Dynamic Programming, 2nd Edition, Athena Scientific, Belmont, MA  (can be  downloaded from the author's website).

\ref[Ber19] Bertsekas, D.\ P., 2019.\ Reinforcement Learning and Optimal Control, Athena Scientific, Belmont, MA.

\ref[Ber20] Bertsekas, D.\ P., 2020.\ Rollout, Policy Iteration, and Distributed Reinforcement Learning, Athena Scientific, Belmont, MA.

\ref [Kle68] Kleinman, D.\ L., 1968.\  ``On an Iterative Technique for Riccati
Equation Computations," IEEE Trans.\ Aut.\  Control, Vol.\ AC-13, pp.\ 114-115.

\ref[Hew71] Hewer, G., 1971.\ ``An Iterative Technique for the Computation of the Steady State Gains for the Discrete Optimal Regulator," IEEE Trans.\ on Automatic Control, Vol.\ 16, pp.\ 382-384. 

\ref[PoA69] Pollatschek, M.\ A., and Avi-Itzhak, B., 1969.\ ``Algorithms for Stochastic Games with Geometrical Interpretation," Management Science,  Vol.\ 15, pp.\ 399-415.

\ref [PuB78] Puterman, M.\ L., and Brumelle, S.\ L., 1978.\  ``The Analytic Theory of
Policy Iteration," in Dynamic Programming and Its Applications, M.\ L.\ Puterman
(ed.), Academic Press, N.\ Y.

\ref [PuB79] Puterman, M.\ L., and Brumelle, S.\ L., 1979.\  ``On the Convergence of Policy Iteration in Stationary Dynamic Programming," Mathematics of Operations Research, Vol.\ 4, pp.\ 60-69.

\ref [Put94] Puterman, M.\ L., 1994.\  Markovian Decision Problems, J.\ Wiley, N.\ Y.

\end